# HARRIS PROCESSES

**\* Sherly Sebastian**

*Department of Statistics, Vimala College, Trichur-680009, Kerala, India.*
(*Tel: 91-0487 -2331888; e-mail: sebastian_sherly@yahoo.com*)

**Jos M. K**

*Department of Statistics, St. Thomas College, Trichur-680001, India.*
(*e-mail: jos_kuriakku@yahoo.com*)

**Sandhya E**

*Department of Statistics, Prajyoti Niketan College, Pudukkad , Trichur-680301, India.*
(*e-mail: esandhya@hotmail.com*)

and

**Raju N**

*Department of Statistics, University of Calicut, Calicut University PO-673635, India.*
(*e-mail: na_raju@yahoo.com*)

**Abstract**

In this paper, we develop two stochastic models where the variable under consideration follows Harris distribution. The mean and variance of the processes are derived and the processes are shown to be non-stationary. In the second model, starting with a Poisson process, an alternate way of obtaining Harris process is introduced.

\* Address for communication

## 1. Introduction

Harris [1] introduced a probability generating function (p.g.f) given by

$$P(s) = \frac{s}{\left(m - (m-1)s^k\right)^{1/k}}, \quad k > 0 \text{ integer}, \quad m > 1.$$

He discussed this p.g.f while considering a simple discrete branching process where a particle either splits into ($k+1$) identical particles or remains the same during a short time interval $\Delta t$. The probability distribution corresponding to the above p.g.f is called Harris distribution. A basic study of this distribution has been done in Sandhya *et al*. ( 2005, available at *http://arxiv.org/abs/math.ST/0506220* ). The probability mass function (p.m.f) of a random variable (r.v) $X$ following a Harris distribution is given by

$$f(x) = \binom{(1/k)+n-1}{n} \left(\frac{1}{m}\right)^{1/k} \left(1 - \frac{1}{m}\right)^n, \quad x = 1 + nk, \quad n = 0,1,2,\ldots$$

where $k > 0$, integer and $m > 1$ are the parameters. We use the notation $X \sim H_1(m,k,1/k)$ to denote the r.v following Harris distribution with parameters $m$ and $k$. It may be noted that unlike the well known standard discrete distributions, Harris distribution is concentrated on $\{1, 1+k, 1+2k, \ldots\}$ where $k$ is fixed.

In this paper, two stochastic models are introduced. In the first model the rate $\lambda$ of the process is considered as a function of $n$, the number of occurrences of the event at the instant. Then the r.v under consideration follows a Harris distribution and the corresponding stochastic process is termed as a Harris process. It is also proved that Yule – Furry process (see Bhat [2], section 8.1.3) is a particular case of Harris process. In the second model, the Poisson rate $\lambda$ follows a gamma distribution and Harris distribution is obtained as a mixture by considering a linear function of Poisson variable.

The corresponding stochastic process is a Harris process. It may be noted that Greenwood and Yule [3] derived the relationship between the Poisson and negative binomial distribution by considering the intensity parameter $\lambda$ of the Poisson process to be a gamma distributed r.v . Also, Barndorff - Nielsen and Yeo [4] defined a negative binomial process as a conditional Poisson process whose intensity function $\lambda(t)$ is a gamma process.

We arrange the paper as follows: Section 2 contains the description and analysis of model 1. Mean and variance of the process are derived. Further, Yule-Furry process is shown to be a particular case of Harris process. In section 3 model 2 is presented. This section provides an alternate approach for obtaining Harris process.

## 2. Model 1

Consider a company producing costly articles like car, television, refrigerator etc. In order to enhance the selling of products the company provides incentives to its marketing executives for their performance. One of the main problems in a manufacturing unit is the uncertainty in the selling volume which depends on the selling rate $\lambda$. Suppose that the company has decided to market its product by taking a policy as follows:

A person is allowed to become a sales executive of the company when he sells one article and he will be given one incentive. Thereafter, for the selling of each $k$ articles, ($k>0$ integer) one additional incentive will be given. Let $I(t)$ denote the number of additional incentives obtained in the interval of duration $t$ starting from an initial epoch $t = 0$. The family of r.v s $\{I(t), t \geq 0\}$ is a stochastic process in continuous time with discrete state space $\{0,1,2,…\}$. Now, let $N(t)$ denote the minimum number of articles sold in order to have $(n+1)$ incentives in the time interval $(0, t]$, where $n = 0,1,2, \ldots$.

The family of r.v s $\{N(t), t \geq 0\}$ is a stochastic process with $N(0) = 1$. Here the time $t$ is continuous and the state space of $N(t)$, $\{1, 1+k, 1+2k, \ldots\}$ is discrete and integer-valued. Thus the whole system is a two-dimensional continuous time stochastic process $\{N(t), I(t), t \geq 0\}$ in the state space $\{(1+nk, n); n \geq 0, k > 0 \text{ integers}\}$.

Let $P_{nk+1}(t)$ be the probability that the r.v $N(t)$ assumes the value $nk+1$.

$$P_{nk+1}(t) = P\{N(t) = nk+1\}, \quad n = 0, 1, 2, \ldots, \quad k > 0 \text{ integer.}$$

$P_{nk+1}(t)$ is a function of the time $t$ and $\sum_{n=0}^{\infty} P_{nk+1}(t) = 1$. $\{P_{nk+1}(t)\}$ represents the probability distribution of the r.v $N(t)$ for every value of $t$.

Let the probability of getting $r$ incentives in $(t, t+\Delta t)$ given that $n$ incentives were obtained by epoch $t$ is given by

$$P(I(\Delta t) = r / I(t) = n) = \lambda_n \Delta t + 0(\Delta t), \quad r = 1$$

$$= 0(\Delta t) \quad r \geq 2$$

$$= 1 - \lambda_n \Delta t + 0(\Delta t) \quad r = 0$$

where $\lambda_n = (nk+1)\lambda$ is a linear function of $n$ and the initial selling rate $\lambda > 0$ is a parameter. Similar models with different linear structures are used in birth and death processes. For examples, see section 6.3. of Ross [5].

**Theorem 1:** Under the above conditions $N(t)$ follows a Harris distribution, i.e. $P_{nk+1}(t)$ is given by the Harris law:

$$P_{nk+1}(t) = \binom{(1/k)+n-1}{n} (e^{-t\lambda k})^{1/k} (1 - e^{-t\lambda k})^n, \quad n = 0, 1, 2, \ldots$$

**Proof:** Consider $P_{nk+1}(t+\Delta t)$ for $n \geq 0$.

For $n \geq 1$, the probability of selling $(nk+1)$ items in order to have $(n+1)$ incentives by time $(t+\Delta t)$ can be written as

$$P_{nk+1}(t+\Delta t) = P_{nk+1}(t)(1-(nk+1)\lambda \Delta t) + P_{(n-1)k+1}(t)((n-1)k+1)\lambda \Delta t + 0(\Delta t) \qquad (1)$$

where as for $n = 0$,

$$P_1(t+\Delta t) = P_1(t)(1-\lambda \Delta t) + 0(\Delta t) \qquad (2)$$

Combining (1) and (2), we get

$$P_{nk+1}(t+\Delta t) = P_{nk+1}(t)(1-(nk+1)\lambda \Delta t) + P_{(n-1)k+1}(t)((n-1)k+1)\lambda \Delta t + 0(\Delta t) \qquad (3)$$

for $n \geq 0$ with $P_{1-k}(t) = 0$. Hence we have

$$\frac{P_{nk+1}(t+\Delta t) - P_{nk+1}(t)}{\Delta t} = ((n-1)k+1)\lambda\, P_{(n-1)k+1}(t) - (nk+1)\lambda\, P_{nk+1}(t) + \frac{0(\Delta t)}{\Delta t}$$

Taking the limit as $\Delta t \to 0$, we get

$$\frac{d}{dt} P_{nk+1}(t) = ((n-1)k+1)\lambda\, P_{(n-1)k+1}(t) - (nk+1)\lambda\, P_{nk+1}(t) \qquad (4)$$

Since the process starts by the selling of the first item $N(0) = 1$. To the system (3), we add the initial condition

$$P_1(0) = 1; \quad P_r(0) = 0 \quad \text{for} \quad r \neq 1. \qquad (5)$$

The probability generating function $P(s, t)$ for the probability distribution $P_{nk+1}(t)$ is

$$P(s,t) = \sum_{n=0}^{\infty} P_{nk+1}(t)\, s^{nk+1}$$

where $P(s, t)$ is a function of $s$ and $t$. Then,

$$\frac{\partial P}{\partial s} = \sum_{n=0}^{\infty} (nk+1) P_{nk+1}(t)\, s^{nk} \quad \text{and} \quad \frac{\partial P}{\partial t} = \sum_{n=0}^{\infty} \frac{\partial}{\partial t} P_{nk+1}(t)\, s^{nk+1}$$

Now, the initial condition (5) can be expressed as

$$P(s,0) = s$$

Multiply (4) by $s^{nk+1}$ and sum over all values of $n$. This gives

$$\frac{\partial P(s,t)}{\partial t} = \lambda s^{k+1} \frac{\partial P(s,t)}{\partial s} - \lambda s \frac{\partial P(s,t)}{\partial s} \ . \quad \text{Thus we have}$$

$$\frac{\partial P(s,t)}{\partial t} + \lambda s (1 - s^k) \frac{\partial P(s,t)}{\partial s} = 0$$

Auxiliary equations are

$$\frac{dt}{1} = \frac{ds}{\lambda s (1 - s^k)} = \frac{dP}{0} \ . \quad \text{For details, see section 1.4 in Amaranath [6].}$$

Solving these equations with the initial conditions $t = 0$ and $P(s,0) = s$ we get

$$P(s,t) = \frac{s}{(e^{t\lambda k} - (e^{t\lambda k} - 1)s^k)^{1/k}} \ , \quad \text{where} \quad e^{t\lambda k} > 1 \tag{6}$$

This is the p.g.f of the Harris distribution $H_1(e^{t\lambda k}, k, 1/k)$. Hence $N(t)$ follows a Harris distribution. Picking out the coefficient of $s^{nk+1}$ on the right side of (6) gives

$$P_{nk+1}(t) = \binom{(1/k)+n-1}{n} (e^{-t\lambda k})^{1/k} (1 - e^{-t\lambda k})^n \ , \quad n = 0, 1, 2, \ldots$$

This completes the theorem.

The mean and variance of this distribution can be obtained from the p.g.f given in (6). Hence we have

**Corollary:** The mean and variance of the process are

$$E\{N(t)\} = e^{t\lambda k} \quad \text{and}$$

$$\mathrm{Var}\{N(t)\} = e^{t\lambda k}(e^{t\lambda k} - 1)k$$

**Remarks:**

(i) The mean number of items sold in order to have $(n+1)$ incentives in an interval of

length $t$ is $e^{t\lambda k}$, so that the mean number of items sold for $(n+1)$ incentives per unit time $(t = 1)$, i.e. in an interval of unit length is $e^{\lambda k}$.

(ii) The mean and variance of $N(t)$ are functions of $t$. Since the distribution of $N(t)$ functionally dependent on $t$, the process $\{ N(t), t \geq 0 \}$ is not stationary.

(iii) We have $N(t) = 1 + k\, I(t)$

$$P(I(t) = n) = P(N(t) = nk + 1)$$

$$= \binom{(1/k)+n-1}{n} (e^{-t\lambda k})^{1/k} (1 - e^{-t\lambda k})^n, \quad n = 0, 1, 2, \ldots$$

which is the p.m.f of the negative binomial distribution $NB(1/k, e^{-t\lambda k})$ defined on $\{0, 1, 2, \ldots\}$. The mean and variance of $I(t)$ are

$E\{I(t)\} = (e^{t\lambda k} - 1) / k$ and

$Var\{I(t)\} = e^{t\lambda k}(e^{t\lambda k} - 1) / k$

$\{N(t), t \geq 0\}$ is a stochastic process which we call 'Harris process'.

**Particular case of Harris process**

When $k = 1$, the distribution of $N(t)$ is given by

$$P_n(t) = P(N(t) = n) = (e^{-t\lambda})(1 - e^{-t\lambda})^{n-1}, \quad n = 1, 2, \ldots$$

which is the p.m.f of the modified geometric or decapitated geometric with parameter $e^{-\lambda t}$ and its p.g.f is given by

$$P(s, t) = \frac{e^{-\lambda t} s}{1 - (1 - e^{-\lambda t})s}$$

The mean and variance of the corresponding process are

$E\{N(t)\} = e^{t\lambda}$ and

$Var\{N(t)\} = e^{t\lambda}(e^{t\lambda k} - 1)$

So when $k = 1$, the stochastic process $\{N(t), t \geq 0\}$ is the Yule-Furry process.

## 3. Model 2

In this section we have developed a stochastic process where the variable under consideration follows a Harris law. We start with a Poisson process with rate $\lambda$ as a gamma distributed r.v. Then a linear function of the Poisson variable leads to Harris process.

Consider a random event $E$ such as the arrival of particles at a counter, occurrence of accidents at a certain place, selling of costly articles from a production unit etc. Let $X(t)$ denote the total number of occurrences of the event $E$ in an interval of duration $t$. Let

$$P^*_n(t) = P(X(t) = n), \qquad n = 0,1,2,\ldots$$

$$\sum_{n=0}^{\infty} P^*_n(t) = 1$$

$\{P^*_n(t)\}$ represents the probability distribution of the r.v $X(t)$ for every value of $t$. The family of r.v s $\{X(t), t \geq 0\}$ is a stochastic process in continuous time $t$. The state space of $X(t)$ is discrete and integer-valued.

Under the postulates for Poisson process (see Karlin and Taylor [7], section 4.1A), it can be proved that $X(t)$ follows Poisson distribution with mean $\lambda t$ where $\lambda > 0$ is a parameter denoting the intensity of selling articles. The probability distribution $\{P^*_n(t)\}$ of $X(t)$ is given by

$$P^*_n(t) = \frac{e^{-\lambda t}(\lambda t)^n}{n!}, \qquad n = 0, 1, 2, \ldots$$

Consider a linear function of $X(t)$, $Z(t) = k X(t) + 1$ where $k > 0$ integer. The family of r.v s $\{Z(t), t \geq 0\}$ is a stochastic process. Here time $t$ is a continuous variable and $Z(t)$ is a non-negative integer-valued r.v taking values 1, 1+k, 1+2k, …

Let $P'_{nk+1}(t) = P(Z(t) = nk+1)$

**Theorem 2:** Under the above conditions $Z(t)$ follows a Harris distribution, provided the mixing distribution is gamma $G(1/k, a)$. i.e. $P'_{nk+1}(t)$ is given by the Harris law:

$$P'_{nk+1}(t) = \binom{(1/k)+n-1}{n} \left(\frac{a}{a+t}\right)^{1/k} \left(\frac{t}{a+t}\right)^n, \quad n = 0, 1, 2, \ldots$$

**Proof:** $P'_{nk+1}(t) = P(Z(t) = nk+1)$

$$= P(X(t) = n)$$

$$= \frac{e^{-\lambda t}(\lambda t)^n}{n!}, \quad n = 0, 1, 2, \ldots \tag{7}$$

When the Poisson rate $\lambda$ is a constant, (7) gives the conditional distribution of $Z(t)$.

Now, assume that $\lambda$ is a gamma r.v with density function

$$f(\lambda) = \frac{a^{1/k}}{\Gamma(1/k)} e^{-a\lambda} \lambda^{(1/k)-1}, \quad k > 0 \text{ integer}, \ a > 0, \ \lambda > 0.$$

Then the mixture is given by

$$P'_{nk+1}(t) = P(Z(t) = nk+1)$$

$$= \int_0^\infty \frac{e^{-\lambda t}(\lambda t)^n}{n!} \frac{a^{1/k}}{\Gamma(1/k)} e^{-a\lambda} \lambda^{(1/k)-1} d\lambda$$

$$= \binom{(1/k)+n-1}{n} \left(\frac{a}{a+t}\right)^{1/k} \left(\frac{t}{a+t}\right)^n$$

$$= \binom{(1/k)+n-1}{n} \left(\frac{1}{m}\right)^{1/k} \left(1-\frac{1}{m}\right)^n, \quad n = 0, 1, 2, \ldots$$

where $k > 0$ is an integer and $m = (a+t)/a > 1$. This is the p.m.f of the Harris distribution with index $1/k$ and mean $(a+t)/a$. i.e. $Z(t) \sim H_1(m,k,1/k)$.

**Corollary:** The mean and variance of the process are

$$E\{Z(t)\} = \frac{a+t}{a} \quad \text{and}$$

$$\text{Var}\{Z(t)\} = \frac{(a+t)tk}{a^2}$$

**Remark:** The mean and variance of $Z(t)$ are functions of $t$. Since the distribution of $Z(t)$ functionally dependent on $t$, the process $\{Z(t), t \geq 0\}$ is non-stationary.